# Stochastic Frank-Wolfe Methods for Nonconvex Optimization

Sashank J. Reddi
sjakkamr@cs.cmu.edu
Carnegie Mellon University

Suvrit Sra
suvrit@mit.edu
Massachusetts Institute of Technology

Barnabás Póczós
bapoczos@cs.cmu.edu
Carnegie Mellon University

Alex Smola
alex@smola.org
Carnegie Mellon University

## Abstract

We study Frank-Wolfe methods for *nonconvex* stochastic and finite-sum optimization problems. Frank-Wolfe methods (in the convex case) have gained tremendous recent interest in machine learning and optimization communities due to their projection-free property and their ability to exploit structured constraints. However, our understanding of these algorithms in the nonconvex setting is fairly limited. In this paper, we propose nonconvex stochastic Frank-Wolfe methods and analyze their convergence properties. For objective functions that decompose into a finite-sum, we leverage ideas from variance reduction techniques for convex optimization to obtain new variance reduced nonconvex Frank-Wolfe methods that have provably faster convergence than the classical Frank-Wolfe method. Finally, we show that the faster convergence rates of our variance reduced methods also translate into improved convergence rates for the stochastic setting.

## 1 Introduction

We study optimization problems of the form:

$$\min_{x \in \Omega} F(x) := \begin{cases} \mathbb{E}_z[f(x,z)], & \text{(stochastic)} \\ \frac{1}{n}\sum_{i=}^n f_i(x), & \text{(finite-sum)}. \end{cases} \quad (1)$$

We assume that $F$, $f$, and $f_i$ ($i \in \{1, \ldots, n\} \triangleq [n]$) are all differentiable, but possibly *nonconvex*; the domain $\Omega$ is *convex* and *compact*.

Problems of this form are at the heart of machine learning and statistics; for instance, the finite-sum problem arises under the name empirical loss minimization and M-estimation. Examples of such problems include multiclass classification, matrix learning, recommendation systems (Jaggi, 2013; Hazan and Kale, 2012; Hazan and Luo, 2016; Harchaoui et al., 2014).

Within convex optimization, problem (1) is relatively well-studied. Two particularly popular approaches for solving it are: (i) Projected stochastic gradient descent (SGD); and (b) the Frank-Wolfe (FW) method. At each iteration, SGD takes a step in a direction opposite to a stochastic approximation of the gradient $\nabla F$ and uses projection onto $\Omega$ to ensure feasibility. While computing a stochastic approximation to $\nabla F$ is usually inexpensive, in many real settings, the cost projecting onto $\Omega$ can be very high (e.g., projecting onto the trace-norm ball, onto base polytopes in submodular minimization (Fujishige and Isotani, 2011)); and in extreme cases projection can even be computationally intractable (Collins et al., 2008).

In such cases, projection based methods like SGD become impractical. This difficulty underlies the recent surge of interest in Frank-Wolfe methods (Frank and Wolfe, 1956; Jaggi, 2013) (also known as conditional gradient), due to their projection-free property. In particular, FW methods avoid the expensive projection operation and require just a *linear oracle* that solves problems of the form $\min_{x \in \Omega} \langle x, g \rangle$ at each iteration.

Despite the remarkable success of FW approaches in the convex setting, including stochastic problems (Hazan and Luo, 2016), their applicability and *non-asymptotic convergence* for *nonconvex* optimization is largely unstudied. Even for SGD, it is only recently that non-asymptotic convergence analysis for nonconvex optimization was obtained (Ghadimi and Lan, 2013; Ghadimi et al., 2014). More recently, Reddi et al. (2016a;b) obtained variance reduced stochastic methods that converge faster than SGD in the *nonconvex* finite-sum setting.

Similar fast variants of FW for nonconvex problems are not known. Given the vast importance of nonconvex models in machine learning (e.g., in deep learning) and the need to incorporate non-trivial constraints in such models, it is imperative to develop scalable, projection-free methods. This paper presents new FW methods towards this goal. Our main contributions are summarized below, while the key complexity results are listed in Figure 1.

**Main Contributions**. For the nonconvex stochastic setting in (1), we propose a stochastic Frank-Wolfe algorithm (SFW), and provide its convergence analysis. For the nonconvex finite-sum setting, we propose two variance reduced (VR) algorithms: SVFW and SAGAFW, based on the popular VR algorithms SVRG and SAGA, respectively. We show that by carefully selecting the parameters of these algorithms, we can attain faster convergence rates than the deterministic FW. In particular, we prove that SVFW and SAGAFW are faster than deterministic FW by a factor of $n^{1/3}$ and $n^{2/3}$ respectively, where $n$ is the number of component functions in the finite-sum (see (1)). Furthermore, leveraging these variance reduced methods, we propose two algorithms, SVFW-S and SAGAFW-S, for the nonconvex stochastic setting, with faster convergence rates than SFW.

To our knowledge, our work presents the first theoretical improvement for stochastic variants of Frank-Wolfe in the context of nonconvex optimization.

## 1.1 Related Work

The classical Frank-Wolfe method (Frank and Wolfe, 1956) using line-search was analyzed for smooth convex functions $F$ and polyhedral domains $\Omega$. Here, a convergence rate of $O(1/\epsilon)$ to ensure $F(x) - F^* \leq \epsilon$ was proved without additional conditions (Frank and Wolfe, 1956; Jaggi, 2013). There have been several recent works on improving the convergence rates under additional assumptions (Garber and Hazan, 2015; Lacoste-Julien and Jaggi, 2015). More recently, Hazan and Luo (2016) proposed stochastic variants of FW for convex problems of form (1), and showed theoretical improvements over the classical Frank-Wolfe method.

The literature on nonconvex Frank-Wolfe is relatively small. The work (Bertsekas, 1995) proves asymptotic convergence of FW to a stationary point; though, no convergence rates are provided. To the best of our knowledge, Yu et al. (2014) is the first to provide convergence rates for FW-type algorithm in the nonconvex setting. Very recently, Lacoste-Julien (2016) provided a (non-asymptotic) convergence rate of $O(1/\epsilon^2)$ for nonconvex FW with adaptive step sizes. However, as we shall see later, implementation of classical FW for (1) is expensive (or impossible in the pure stochastic case) since it requires calculation of the gradient $\nabla F$ at each iteration. We show that our stochastic variants are provably faster than the existing FW methods.

In the nonconvex setting, most of the work on stochastic methods focuses on SGD (Ghadimi and Lan, 2013; Ghadimi et al., 2014) and analyzes convergence to stationary points. For the finite-sum setting, we build on recent variance reduction techniques (Johnson and Zhang, 2013; Defazio et al., 2014; Schmidt et al., 2013), which were first proposed for solving unconstrained convex problems of form (1). Projected variants to handle constraints were studied in (Defazio et al., 2014; Xiao and Zhang, 2014). More recently, Reddi et al. (2016a;b;c) provided nonconvex variants of these methods

| Algorithm | SFO/IFO Complexity | LO Complexity |
|-----------|-------------------|---------------|
| Frank-Wolfe | $O(n/\epsilon^2)$ | $O\left(1/\epsilon^2\right)$ |
| SFW | $O\left(1/\epsilon^4\right)$ | $O\left(1/\epsilon^2\right)$ |
| SVFW | $O(n + n^{2/3}/\epsilon^2)$ | $O(1/\epsilon^2)$ |
| SAGAFW | $O(n + n^{1/3}/\epsilon^2)$ | $O(1/\epsilon^2)$ |
| SVFW-S | $O(1/\epsilon^{10/3})$ | $O(1/\epsilon^2)$ |
| SAGAFW-S | $O(1/\epsilon^{8/3})$ | $O(1/\epsilon^2)$ |

Figure 1: Table comparing the *best* SFO/IFO and LO complexity of algorithms discussed in the paper (for the nonconvex setting). Here, SFW, SVFW-S and SAGAFW-S are algorithms for the stochastic setting, while FW, SVFW and SAGAFW are algorithms for the finite-sum setting. The complexity is measured by the number of oracle calls required to achieve an $\epsilon$-accurate solution (see Section 2 for definitions of SFO/IFO and LO complexity). The complexity of FW is from (Lacoste-Julien, 2016). The results marked in red are contributions of this paper. For clarity, we hide the dependence of SFO/IFO and LO complexity on the initial point and few parameters related to the function $F$ and domain $\Omega$.

that converge provably faster than both SGD and its deterministic counterpart.

## 2 Preliminaries

As stated above, we study two different problem settings: (1) *stochastic*, where $F(x) = \mathbb{E}_z[f(x,z)]$ and $z$ is random variable whose distribution $\mathcal{P}$ is supported on $\Xi \subset \mathbb{R}^p$; and (2) *finite-sums*, where $F(x) = \frac{1}{n}\sum_{i=1}^n f_i(x)$.

For the stochastic setting, we assume that $F$ is *L-smooth*, i.e., its gradient is Lipschitz continuous with constant $L$, so

$$\|\nabla F(x) - \nabla F(y)\| \leq L\|x - y\|, \quad \forall\ x, y \in \Omega.$$

Here $\|.\|$ denotes the $\ell^2$-norm. Furthermore, for the stochastic setting, we also assume the function $f$ is $G$-Lipschitz i.e., $\|\nabla f(x, z)\| \leq G$ for all $x \in \Omega$ and $z \in \Xi$. Such an assumption is common in the stochastic setting (Ghadimi and Lan, 2013; Hazan and Luo, 2016).

For the finite-sum setting, we assume that the individual functions $f_i$ ($i \in [n]$) are $L$-smooth i.e.,

$$\|\nabla f_i(x) - \nabla f_i(y)\| \leq L\|x - y\|, \quad \forall\ x, y \in \Omega.$$

Note that this implies that the function $F$ is also $L$-smooth. The domain $\Omega \in \mathbb{R}^d$ is assumed to be convex and compact with diameter $D$; i.e., $\|x - y\| \leq D$ for all $x, y \in \Omega$. Such an assumption is common to all Frank-Wolfe methods.

**Convergence criteria.** The criterion used for the convergence analysis is important in nonconvex optimization. For unconstrained problems, the gradient norm $\|\nabla F\|$ is typically used to measure convergence, because $\|\nabla F\| \to 0$ translates into convergence to a stationary point. However, this criterion cannot be used for constrained problems of the form (1). Instead, we use the following quantity, typically referred to as *Frank-Wolfe gap*:

$$\mathcal{G}(x) = \max_{v \in \Omega} \langle v - x, -\nabla F(x) \rangle. \tag{2}$$

For convex functions, the Fw gap provides an upper bound on the suboptimality. For nonconvex functions, the gap $\mathcal{G}(x) = 0$ if and only if $x$ is a stationary point. To state our convergence results we will also need the following bound:

$$\beta \geq \frac{2(F(x_0) - F(x^*))}{LD^2},$$

given some (unspecified) initial point $x_0 \in \Omega$.

**Oracle model.** To compare convergence speed of different algorithms, we use the following black-box oracles:

- *Stochastic First-Order Oracle* (SFO): For a function $F(\cdot) = \mathbb{E}_z[f(.,z)]$ where $z \sim \mathcal{P}$, an SFO takes a point $x$ and returns the pair $(f(x, z'), \nabla f(x, z'))$ where $z'$ is a sample drawn i.i.d. from $\mathcal{P}$ (Nemirovski and Yudin, 1983).

- *Incremental First-Order Oracle* (IFO): For a function $F(\cdot) = \frac{1}{n} \sum_i f_i(.)$, an IFO takes an index $i \in [n]$ and a point $x \in \mathbb{R}^d$, and returns the pair $(f_i(x), \nabla f_i(x))$ (Agarwal and Bottou, 2014).

- *Linear Optimization Oracle* (LO): For a set $\Omega$, an LO takes a direction $d$ and returns $\arg\max_{v \in \Omega} \langle v, d \rangle$.

Throughout the paper, by SFO, IFO and LO complexity of an algorithm, we mean the total number of SFO, IFO and LO calls made by the algorithm to obtain an $\epsilon$-*accurate solution*, i.e., a solution for which $\mathbb{E}[\mathcal{G}(x)] \leq \epsilon$; the expectation is over any randomization as part of the algorithm. For clarity of presentation, we hide the dependence of these complexities on the initial point $F(x_0) - F(x^*)$, Lipschitz constant $G$, and the smoothness constant $L$; we report the dependence on $n$ to highlight its importance.

**Classical Fw.** To place our results in perspective, we begin by recalling the classical Frank-Wolfe (Fw) algorithm (Frank and Wolfe, 1956). Pseudocode for this is presented in Algorithm 1.

---

**Algorithm 1:** Fw $\left(x_0, T, \{\gamma_i\}_{i=0}^{T-1}\right)$

1: **Input:** $x_0 \in \Omega$, number of iterations $T$, $\{\gamma_i\}_{i=0}^{T-1}$ where $\gamma_i \in [0, 1]$ for all $i \in \{0, \ldots, T-1\}$
2: **for** $t = 0$ **to** $T - 1$ **do**
3:     Compute $v_t = \arg\max_{v \in \Omega} \langle v, -\nabla F(x_t) \rangle$
4:     Compute update direction $d_t = v_t - x_t$
5:     $x_{t+1} = x_t + \gamma_t d_t$
6: **end for**

---

Each iteration of Fw entails calculation of the gradient $\nabla F$ and moving towards a minimizer of a linearized objective. Notice that calculation of $\nabla F$ may not be possible in the stochastic setting of (1). Furthermore, even in the finite-sum setting, computing $\nabla F$ requires $n$ IFO calls, rendering the approach useless in large-scale problems, where $n$ is large. For the nonconvex finite-sum setting, the following key result was proved recently (Lacoste-Julien, 2016).

**Theorem 1** (Lacoste-Julien, 2016). *Under appropriate selection of step sizes $\gamma_t$, the IFO and LO complexity of Algorithm 1 to achieve an $\epsilon$-accurate solution in the finite-sum setting are $O(n/\epsilon^2)$ and $O(1/\epsilon^2)$ respectively.*

The key aspect of Theorem 1 is the dependence of IFO complexity on $n$. In particular, when $n$ is large, the IFO complexity $O(n/\epsilon^2)$ shown by the theorem becomes prohibitively expensive; thus, undermining the benefits of Fw over competitors like projected SGD. In the next section, we tackle this drawback and develop faster nonconvex stochastic and variance reduced Fw methods.

**Algorithm 2:** Nonconvex SFW $(x_0, T, \{\gamma_i\}_{i=0}^{T-1}, \{b_i\}_{i=0}^{T-1})$

1: **Input:** $x_0 \in \Omega$, number of iterations $T$, $\{\gamma_i\}_{i=0}^{T-1}$ where $\gamma_i \in [0,1]$ for all $i \in \{0, \ldots, T-1\}$, minibatch size $\{b_i\}_{i=0}^{T-1}$
2: **for** $t = 0$ to $T-1$ **do**
3:     Uniformly randomly pick i.i.d samples $\{z_1^t, \ldots, z_{b_t}^t\}$ according to the distribution $\mathcal{P}$.
4:     Compute $v_t = \arg\max_{v \in \Omega} \langle v, -\frac{1}{b_t} \sum_{i=1}^{b_t} \nabla f(x_t, z_i) \rangle$
5:     Compute update direction $d_t = v_t - x_t$
6:     $x_{t+1} = x_t + \gamma_t d_t$
7: **end for**
8: **Output:** Iterate $x_a$ chosen uniformly random from $\{x_t\}_{t=0}^{T-1}$.

## 3 Algorithms

In this section, we describe FW algorithms for solving (1). In particular, we explore stochastic and variance reduced versions of the classical FW method, for the stochastic and finite-sum settings, respectively. We defer the discussion on comparison of the convergence rates to Section 5.

### 3.1 Stochastic Setting

We first investigate the convergence of FW in the stochastic setting. As mentioned earlier, the classical FW method (Algorithm 1) requires calculation of the full gradient $\nabla F(x)$, which is typically impossible to compute in the stochastic setting. For convex problems, Hazan and Luo (2016) tackle this issue by using the popular Robbins-Monro approximation (Robbins and Monro, 1951) to the gradient. We use a variant of the algorithm for our nonconvex stochastic setting, which we call SFW.

The pseudocode of SFW is listed in Algorithm 2. Note that the samples $\{z_i\}$ are chosen independently according to the distribution $\mathcal{P}$. Thus, $\mathbb{E}_{z_i}[\nabla f(x, z_i)] = \nabla F(x)$, i.e., we obtain an unbiased estimate of the gradient. Also, note that the output in Algorithm 2 is randomly selected from all the iterates of the algorithm. The key parameters of SFW are the step sizes $\{\gamma_i\}_{i=0}^{T-1}$ and the minibatch sizes $\{b_t\}$. These parameters must be chosen appropriately in order to ensure convergence of the algorithm (see Theorem 2). For our analysis, we assume that the function $f$ is $G$-Lipschitz i.e., we have $\max_{x \in \Omega, z \in \Xi} \|\nabla f(x, z)\| \leq G$. This bound on the gradient is crucial for our convergence analysis.

We prove the following key result for nonconvex SFW.

**Theorem 2.** *Consider the stochastic setting of* (1) *where $f$ is $G$-Lipschitz and $F$ is $L$-smooth. Then, the output $x_a$ of Algorithm 2 with parameters $\gamma_t = \gamma = \sqrt{\frac{2(F(x_0) - F(x^*))}{TLD^2\beta}}$, $b_t = b = T$ for all $t \in \{0, \ldots, T-1\}$, satisfies the following bound:*

$$\mathbb{E}[\mathcal{G}(x_a)] \leq \frac{D}{\sqrt{T}} \left( G + \sqrt{\frac{2L(F(x_0) - F(x^*))}{\beta}}(1 + \beta) \right),$$

*where $x^*$ is an optimal solution to (stochastic) problem* (1).

*Proof.* First observe the following upper bound:

$$\begin{aligned} F(x_{t+1}) &\leq F(x_t) + \langle \nabla F(x_t), x_{t+1} - x_t \rangle + \frac{L}{2} \|x_{t+1} - x_t\|^2 \\ &= F(x_t) + \langle \nabla F(x_t), \gamma(v_t - x_t) \rangle + \frac{L}{2} \|\gamma(v_t - x_t)\|^2 \\ &\leq F(x_t) + \gamma \langle \nabla F(x_t), v_t - x_t \rangle + \frac{LD^2 \gamma^2}{2}. \end{aligned} \qquad (3)$$

The first inequality follows since $F$ is L-smooth (see Lemma 1). The equality is due to the fact that $x_{t+1} - x_t = \gamma(v_t - x_t)$. The second inequality holds because $v_t, x_t \in \Omega$ and because the diameter of $\Omega$ is $D$.

Next, we introduce the following quantity:

$$\hat{v}_t := \arg\max_{v \in \Omega} \langle v, -\nabla F(x_t) \rangle, \tag{4}$$

which is used purely for our analysis and is *not* part of the algorithm. For brevity, we use $\nabla_t$ to denote $\frac{1}{b}\sum_{i=1}^{b} f(x_t, z_i^t)$.

Rewriting inequality (3) using this quantity, we see that

$$\begin{aligned}
F(x_{t+1}) &\leq F(x_t) + \gamma \langle \nabla_t, v_t - x_t \rangle \\
&\quad + \gamma \langle \nabla F(x_t) - \nabla_t, v_t - x_t \rangle + \frac{LD^2\gamma^2}{2} \\
&\leq F(x_t) + \gamma \langle \nabla_t, \hat{v}_t - x_t \rangle \\
&\quad + \gamma \langle \nabla F(x_t) - \nabla_t, v_t - x_t \rangle + \frac{LD^2\gamma^2}{2} \\
&= F(x_t) + \gamma \langle \nabla F(x_t), \hat{v}_t - x_t \rangle \\
&\quad + \gamma \langle \nabla F(x_t) - \nabla_t, v_t - \hat{v}_t \rangle + \frac{LD^2\gamma^2}{2} \\
&= F(x_t) - \gamma \mathcal{G}(x_t) + \gamma \langle \nabla F(x_t) - \nabla_t, v_t - \hat{v}_t \rangle + \frac{LD^2\gamma^2}{2} \\
&\leq F(x_t) - \gamma \mathcal{G}(x_t) + D\gamma \|\nabla F(x_t) - \nabla_t\| + \frac{LD^2\gamma^2}{2}. \tag{5}
\end{aligned}$$

The second inequality follows from the optimality of $v_t$ in Algorithm 2, while the third inequality follows from recalling that $\mathcal{G}(x_t) = \max_{v \in \Omega} \langle v - x_t, -\nabla F(x_t) \rangle = \langle \hat{v}_t - x_t, -\nabla F(x_t) \rangle$, which holds due to the optimality of $\hat{v}_t$ in (4). The last inequality follows from Cauchy-Schwarz and the fact that the diameter of the feasible set $\Omega$ is bounded by $D$.

Taking expectations and using Lemma 2 in (5) we obtain the following important bound:

$$\mathbb{E}[F(x_{t+1})] \leq \mathbb{E}[F(x_t)] - \gamma \mathbb{E}[\mathcal{G}(x_t)] + \frac{GD\gamma}{\sqrt{b}} + \frac{LD^2\gamma^2}{2}.$$

Summing over $t$ and telescoping, we then obtain the upper-bound

$$\begin{aligned}
\gamma \sum_{t=0}^{T-1} \mathbb{E}[\mathcal{G}(x_t)] &\leq F(x_0) - \mathbb{E}[F(x_T)] + \frac{TGD\gamma}{\sqrt{b}} + \frac{TLD^2\gamma^2}{2} \\
&\leq F(x_0) - F(x^*) + \frac{TGD\gamma}{\sqrt{b}} + \frac{TLD^2\gamma^2}{2}.
\end{aligned}$$

The latter inequality follows from the optimality of $x^*$. Using the definition of the output $x_a$ of Algorithm 2 and the parameters specified in the theorem statement, we get

$$\begin{aligned}
\mathbb{E}[\mathcal{G}(x_a)] &\leq \frac{F(x_0) - F(x^*)}{T\gamma} + \frac{GD}{\sqrt{b}} + \frac{LD^2\gamma}{2} \\
&\leq \frac{D}{\sqrt{T}}\left(G + \sqrt{\frac{2L(F(x_0) - F(x^*))}{\beta}}(1 + \beta)\right), \tag{6}
\end{aligned}$$

which concludes the proof of the theorem. $\square$

An immediate consequence of Theorem 2 is the following complexity result for SFW.

**Corollary 1.** *Under the setting of Theorem 2, the SFO complexity and LO complexity of Algorithm 2 are $O(1/\epsilon^4)$ and $O(1/\epsilon^2)$, respectively.*

*Proof.* The proof follows upon observing that $O(1/\epsilon^2)$ minibatch size is required at each iteration of the algorithm, and noting that as per Theorem 2 $O(1/\epsilon^2)$ iterations are required to achieve an $\epsilon$-accurate solution. □

Note that the SFO and LO complexity of nonconvex SFW is similar to that of online FW (Hazan and Kale, 2012) and slightly worse than complexity of SFW for convex problems (Hazan and Luo, 2016). Furthermore, for simplicity of analysis, we used a fixed step size and minibatch size. One can derive an essentially similar result using a decreasing step size and increasing minibatch size.

It is important to emphasize that the above results also apply to the finite-sum setting. In particular, when the distribution $\mathcal{P}$ is the empirical measure, then the convergence result in Theorem 2 also provides convergence rates for the finite-sum case. However, as we will see shortly, these convergence rates can be improved significantly by using variance reduction techniques.

### 3.2 Finite-sum Setting

In this section, we consider the finite-sum setting of (1). We show that by building on ideas from variance reduction for SGD, one can significantly improve the convergence rates. The key idea is to use a variance reduced approximation of the gradient (Johnson and Zhang, 2013; Defazio et al., 2014). We analyze two different algorithms for the finite-sum setting. Our first algorithm (SVFW) is based on the convex method of (Hazan and Luo, 2016) adapted to the nonconvex case. Our second algorithm (SAGAFW) is based on another variance reduction technique called SAGA (Defazio et al., 2014).

## SVFW Algorithm

Pseudocode of our first method (SVFW) is presented in Algorithm 3. Similar to (Johnson and Zhang, 2013) and (Hazan and Luo, 2016), nonconvex SVFW is also epoch-based. At the end of each epoch, the full gradient is computed at the current iterate. This gradient is used for controlling the variance of the stochastic gradients in the inner loop. For epoch size $m = 1$, SVFW reduces to the classic FW algorithm. In general, the epoch size $m$ is chosen such that the total number of IFO calls per epoch is $\Theta(n)$. This ensures that the cost of computing the full gradient at the end of each epoch is amortized. To enable a fair comparison with SFW, we assume that the total number of inner iterations across all epochs in Algorithm 3 is $T$.

We prove the following key result for Algorithm 3. For ease of exposition, we assume that the total number of inner iterations $T$ is a multiple of $m$.

**Theorem 3.** *Consider the finite-sum setting of (1) where the functions $\{f_i\}_{i=1}^n$ are L-smooth. Then, the output $x_a$ of Algorithm 3 with parameters $\gamma_t = \gamma = \sqrt{\frac{F(x_0)-F(x^*)}{TLD^2\beta}}$ and $b_t = b = m^2$ for all $t \in \{0, \ldots, m-1\}$, satisfies*

$$\mathbb{E}[\mathcal{G}(x_a)] \leq \frac{2D}{\sqrt{T\beta}}\sqrt{L(F(x_0) - F(x^*))}(1 + \beta),$$

*where $x^*$ is an optimal solution of (1) and $x_a$ is the output of Algorithm 4.*

*Proof.* We first analyze the convergence properties of iterates within an epoch. Suppose that the current epoch is $s + 1$. For brevity, we drop the symbol $s$ from $x_t^{s+1}$, $\tilde{x}^s$ and $\tilde{g}^s$, whenever it safe to do so given the context. The first part of the proof is similar to that of Theorem 2. For the sake of completeness, we provide the details here. We again use the quantity $\hat{v}_t = \arg\max_{v \in \Omega}\langle v, -\nabla F(x_t)\rangle$, as before, purely for our analysis.

**Algorithm 3:** SVFW $\left(x_0, T, m, \{\gamma_i\}_{i=0}^{m-1}, \{b_i\}_{i=0}^{m-1}\right)$

1: **Input:** $x_m^0 = x_0 \in \Omega$, epoch size $m$, number of epochs $S = \lceil T/m \rceil$, $\{\gamma_i\}_{i=0}^{m-1}$ where $\gamma_i \in [0,1]$ for all $i \in \{0, \ldots, m-1\}$, minibatch size $\{b_i\}_{i=0}^{m-1}$
2: **for** $s = 0$ to $S - 1$ **do**
3:     Let $\tilde{x}^s = x_m^s$
4:     Compute $\tilde{g}^s = \nabla F(\tilde{x}^s) = \frac{1}{n}\sum_{i=1}^n f(\tilde{x}^s)$
5:     **for** $t = 0$ to $m - 1$ **do**
6:         Uniformly randomly (with replacement) select subset $I_t = \{i_1, \ldots, i_{b_t}\}$ from $[n]$.
7:         Compute $v_t^{s+1} = \arg\max_{v \in \Omega} \langle v, -\frac{1}{b_t}(\sum_{i \in I_t} \nabla f_i(x_t^{s+1}) - f_i(\tilde{x}^s) + \tilde{g}^s)\rangle$
8:         Compute update direction $d_t^{s+1} = v_t^{s+1} - x_t^{s+1}$
9:         $x_{t+1}^{s+1} = x_t^{s+1} + \gamma_t d_t^{s+1}$
10:     **end for**
11: **end for**
12: **Output:** Iterate $x_a$ chosen uniformly random from $\{\{x_t^{s+1}\}_{t=0}^{m-1}\}_{s=0}^{S-1}$.

---

For the $t^{th}$ iteration within the epoch $s$, we have

$$F(x_{t+1}) \leq F(x_t) + \langle \nabla F(x_t), \gamma(v_t - x_t)\rangle + \frac{L}{2}\|\gamma(v_t - x_t)\|^2$$
$$\leq F(x_t) + \gamma \langle \nabla F(x_t), v_t - x_t\rangle + \frac{LD^2\gamma^2}{2}. \tag{7}$$

This is due to Lemma 1 and definition of $x_{t+1}$ in Algorithm 3. For brevity, we use $\tilde{\nabla}_t$ to denote $\frac{1}{b_t}(\sum_{i \in I_t} \nabla f_i(x_t) - f_i(\tilde{x}) + \tilde{g})$. Rewriting, we then obtain

$$F(x_{t+1}) \leq F(x_t) + \gamma \langle \tilde{\nabla}_t, v_t - x_t\rangle$$
$$+ \gamma \langle \nabla F(x_t) - \tilde{\nabla}_t, v_t - x_t\rangle + \frac{LD^2\gamma^2}{2}$$
$$\leq F(x_t) + \gamma \langle \tilde{\nabla}_t, \hat{v}_t - x_t\rangle$$
$$+ \gamma \langle \nabla F(x_t) - \tilde{\nabla}_t, v_t - x_t\rangle + \frac{LD^2\gamma^2}{2}$$
$$\leq F(x_t) + \gamma \langle \nabla F(x_t), \hat{v}_t - x_t\rangle$$
$$+ \gamma \langle \nabla F(x_t) - \tilde{\nabla}_t, v_t - \hat{v}_t\rangle + \frac{LD^2\gamma^2}{2}$$
$$\leq F(x_t) - \gamma \mathcal{G}(x_t) + D\gamma\|\nabla F(x_t) - \tilde{\nabla}_t\| + \frac{LD^2\gamma^2}{2}. \tag{8}$$

The second inequality is due to the optimality of $v_t$ in Algorithm 3. The last inequality is due to the definition of $\mathcal{G}(x_t)$, the diameter of set $\Omega$, and an application of Cauchy-Schwarz inequality. Note that the above inequality is similar to (5), except for the crucial difference in the term $\nabla F(x_t) - \tilde{\nabla}_t$ (instead of $\nabla F(x_t) - \nabla_t$ in (5)). As we shall see shortly, this term has much lower variance, which ultimately leads to faster convergence rates.

Taking expectations and using Lemma 3 in inequality (8) we obtain the bound

$$\mathbb{E}[F(x_{t+1})] \leq \mathbb{E}[F(x_t)] - \gamma\mathbb{E}[\mathcal{G}(x_t)]$$
$$+ \frac{LD\gamma}{\sqrt{b}}\mathbb{E}[\|x_t - \tilde{x}\|] + \frac{LD^2\gamma^2}{2}. \tag{9}$$

To aid further analysis, we introduce the following Lyapunov function:

$$R_t = \mathbb{E}[F(x_t) + c_t\|x_t - \tilde{x}\|],$$

where $c_m = 0$ and $c_t = c_{t+1} + (LD\gamma)/\sqrt{b}$ for all $t \in \{0, \ldots, m-1\}$. Using the relationship in (9), we see that

$$\begin{aligned}
R_{t+1} &= \mathbb{E}[F(x_{t+1}) + c_{t+1}\|x_{t+1} - \tilde{x}\|] \\
&\leq \mathbb{E}[F(x_t)] - \gamma \mathbb{E}[\mathcal{G}(x_t)] + \frac{LD\gamma}{\sqrt{b}}\mathbb{E}[\|x_t - \tilde{x}\|] \\
&\quad + \frac{LD^2\gamma^2}{2} + c_{t+1}\mathbb{E}[\|x_{t+1} - \tilde{x}\|] \\
&\leq \mathbb{E}[F(x_t)] - \gamma \mathbb{E}[\mathcal{G}(x_t)] + \frac{LD\gamma}{\sqrt{b}}\mathbb{E}[\|x_t - \tilde{x}\|] \\
&\quad + \frac{LD^2\gamma^2}{2} + c_{t+1}\mathbb{E}[\|x_{t+1} - x_t\| + \|x_t - \tilde{x}\|] \\
&\leq R_t - \gamma \mathbb{E}[\mathcal{G}(x_t)] + \frac{LD^2\gamma^2}{2} + c_{t+1}D\gamma.
\end{aligned} \qquad (10)$$

The second inequality follows from the triangle inequality, while the last inequality holds because: (a) $c_t = c_{t+1} + (LD\gamma)/\sqrt{b}$, and (b) $\|x_{t+1} - x_t\| = \gamma\|v_t - x_t\| \leq D\gamma$ (recall the definition of diameter of $\Omega$). Telescoping over all the iterations within an epoch, we obtain

$$\begin{aligned}
R_m &\leq R_0 - \gamma \sum_{t=0}^{m-1} \mathbb{E}[\mathcal{G}(x_t)] + \frac{LmD^2\gamma^2}{2} + D\gamma \sum_{t=1}^{m} c_t \\
&= R_0 - \gamma \sum_{t=0}^{m-1} \mathbb{E}[\mathcal{G}(x_t)] + \frac{LmD^2\gamma^2}{2} + \frac{L(m-1)mD^2\gamma^2}{2\sqrt{b}}.
\end{aligned} \qquad (11)$$

The equality follows from the relationship $c_t = c_{t+1} + (LD\gamma)/\sqrt{b}$. Since $c_m = 0$ and $x_0^{s+1} = \tilde{x}^s = x_m^s$ (in Algorithm 3), from (11) we obtain

$$\begin{aligned}
\mathbb{E}[F(x_m^{s+1})] &\leq \mathbb{E}[F(x_m^{s+1})] - \gamma \sum_{t=0}^{m-1} \mathbb{E}[\mathcal{G}(x_t^{s+1})] \\
&\quad + \frac{LmD^2\gamma^2}{2} + \frac{L(m-1)mD^2\gamma^2}{2\sqrt{b}}.
\end{aligned}$$

Now telescoping over all epochs, we obtain

$$\begin{aligned}
\mathbb{E}[F(x_m^S)] &\leq F(x_0) - \gamma \sum_{s=0}^{S-1}\sum_{t=0}^{m-1} \mathbb{E}[\mathcal{G}(x_t^{s+1})] \\
&\quad + \frac{TLD^2\gamma^2}{2} + \frac{TL(m-1)D^2\gamma^2}{2\sqrt{b}}.
\end{aligned}$$

Rearranging this inequality and using the definition of the output in Algorithm 3, we finally obtain

$$\begin{aligned}
\mathbb{E}[\mathcal{G}(x_a)] &\leq \frac{F(x_0) - \mathbb{E}[F(x_m^S)]}{T\gamma} + \frac{LD^2\gamma}{2} + \frac{L(m-1)D^2\gamma}{2\sqrt{b}} \\
&\leq \frac{F(x_0) - F(x^*)}{T\gamma} + LD^2\gamma \\
&\leq 2\sqrt{\frac{LD^2(F(x_0) - F(x^*))}{T\beta}}(1+\beta).
\end{aligned}$$

The second inequality follows from the optimality of $x^*$ and because $b = m^2$. The last inequality follows from the choice of $\gamma$ stated in the theorem. This concludes the proof. $\square$

**Algorithm 4:** SAGAFW $(x_0, T, \{\gamma_i\}_{i=0}^{T-1}, \{b_i\}_{i=0}^{T-1})$

1: **Input:** $\alpha_0^i = x_0 \in \Omega$ for all $i \in [n]$, number of iterations $T$, $\{\gamma_i\}_{i=0}^{T-1}$ where $\gamma_i \in [0,1]$ for all $i \in \{0, \ldots, T-1\}$, minibatch size $\{b_i\}_{i=0}^{T-1}$
2: Compute $g_0 = \frac{1}{n} \sum_{i=1}^n \nabla F(\alpha_0^i)$
3: **for** $t = 0$ **to** $T - 1$ **do**
4:     Uniformly randomly (with replacement) select subsets $I_t, J_t$ from $[n]$ of size $b_t$.
5:     Compute $v_t = \arg\max_{v \in \Omega} \langle v, -\frac{1}{b_t}(\sum_{i \in I_t} \nabla f_i(x_t) - f_i(\alpha_t^i) + g_t)\rangle$
6:     Compute update direction $d_t = v_t - x_t$
7:     $x_{t+1} = x_t + \gamma_t d_t$
8:     $\alpha_{t+1}^j = x_t$ for $j \in J_t$ and $\alpha_{t+1}^j = \alpha_t^j$ for $j \notin J_t$
9:     $g_{t+1} = g_t - \frac{1}{n} \sum_{j \in J_t} (\nabla f_j(\alpha_t^j) - \nabla f_j(\alpha_{t+1}^j))$
10: **end for**
11: **Output:** Iterate $x_a$ chosen uniformly random from $\{x_t\}_{t=0}^{T-1}$.

The analysis suggests that the value of $m$ should be set appropriately in Theorem 3 to obtain good convergence rates. If $m$ is small, the IFO complexity of Algorithm 3 is dominated by the step involving calculation of the full gradient at the end of each epoch. On the other hand, if $m$ is large, a large minibatch is used in each step of the algorithm (since $b = m^2$), which increases the IFO complexity. With this intuition, we present following important corollary.

**Corollary 2.** *Under the setting of Theorem 3 and with $m = \lceil n^{1/3} \rceil$, the IFO complexity and LO complexity of Algorithm 2 are $O(n + n^{2/3}/\epsilon^2)$ and $O(1/\epsilon^2)$, respectively.*

*Proof.* We first observe that the total number of IFO calls for an epoch (including those required for calculating the full gradient) is $\Theta(m^3 + n)$. Since $m = \lceil n^{1/3} \rceil$, the total amortized IFO complexity of one iteration within an epoch is $O(m^2) = O(n^{2/3})$. Therefore, the IFO complexity is $O(n + n^{2/3}/\epsilon^2)$. Further, since each inner iteration requires $O(1)$ LO calls, the LO complexity is $O(1/\epsilon^2)$. □

## SAGAFW Algorithm

SVFW is a semi-stochastic algorithm since it requires calculation of the full gradient at the end of each epoch. Below we propose a purely incremental method (SAGAFW) based on the SAGA algorithm of (Defazio et al., 2014). The pseudocode for SAGAFW is presented in Algorithm 4.

A key feature of SAGAFW is that it entirely avoids calculation of full gradients. Instead, it updates the average gradient vector $g_t$ at each iteration. This update requires maintaining additional vectors $\alpha^i$ ($i \in [n]$), and in the worst case such a strategy incurs additional storage cost of $O(nd)$. However, this cost can be reduced to $O(n)$ in several practical cases (refer to (Defazio et al., 2014; Reddi et al., 2016b)).

For SAGAFW, we prove the following key result.

**Theorem 4.** *Consider the finite-sum setting of (1) where functions $\{f_i\}_{i=1}^n$ are L-smooth. Define $\theta(b, n, T) = 1/2 + (2n^{3/2}/Tb^{3/2})$. Then the output $x_a$ of Algorithm 4 with parameters $\gamma_t = \gamma = \sqrt{\frac{F(x_0) - F(x^*)}{TLD^2\theta(b,n,T)\beta}}$ and $b_t = b \leq n$ for all $t \in \{0, \ldots, T-1\}$, satisfies the following:*

$$\mathbb{E}[\mathcal{G}(x_a)] \leq \frac{2D}{\sqrt{T\beta}} \sqrt{L\theta(b,n,T)(F(x_0) - F(x^*))}(1 + \beta),$$

*where $x^*$ is an optimal solution of problem (1) and $x_a$ is the output of Algorithm 4.*

*Proof.* We use the following quantities in our analysis:
$$\check{\nabla}_t = \frac{1}{b_t} \sum_{i \in I_t} \left( \nabla f_i(x_t) - f_i(\alpha_t^i) + g_t \right)$$
$$\hat{v}_t = \arg\max_{v \in \Omega} \langle v, -\nabla F(x_t) \rangle.$$

The first part of our proof is similar to that of Theorem 3. Using essentially the same argument until (8), we have

$$\mathbb{E}[F(x_{t+1})]$$
$$\leq F(x_t) - \gamma \mathcal{G}(x_t) + D\gamma \|\nabla F(x_t) - \check{\nabla}_t\| + \frac{LD^2\gamma^2}{2}$$
$$\leq F(x_t) - \gamma \mathcal{G}(x_t) + \frac{LD\gamma\sqrt{n}}{\sqrt{b}} \frac{1}{n} \sum_{i=1}^n \mathbb{E}\|x_t - \alpha_t^i\| + \frac{LD^2\gamma^2}{2}. \qquad (12)$$

The second inequality is due to Lemma 4. Next, we define the following Lyapunov function:

$$R_t = \mathbb{E}[F(x_t)] + c_t \frac{1}{n} \sum_{i=1}^n \mathbb{E}\|x_t - \alpha_t^i\|,$$

where $c_T = 0$ and $c_t = (1-\rho)c_{t+1} + (LD\gamma\sqrt{n})/\sqrt{b}$ for all $t \in \{0, \ldots, T-1\}$, where $\rho$ is the probability $1 - (1 - 1/n)^b$ of an index $i$ being in $J_t$. We can bound $\rho$ from below as

$$\rho = 1 - \left(1 - \frac{1}{n}\right)^b \geq 1 - \frac{1}{1+(b/n)} = \frac{b/n}{1+b/n} \geq \frac{b}{2n}, \qquad (13)$$

where the first inequality follows from $(1-y)^r \leq 1/(1+ry)$ (which holds for $y \in [0,1]$ and $r \geq 1$), while the second inequality holds because $b \leq n$. Now observe the following:

$$\frac{1}{n} \sum_{i=1}^n \mathbb{E}\|x_{t+1} - \alpha_{t+1}^i\|$$
$$= \frac{1}{n} \sum_{i=1}^n \mathbb{E}\left[\rho\|x_{t+1} - x_t\| + (1-\rho)\|x_{t+1} - \alpha_t^i\|\right]$$
$$\leq \frac{1}{n} \sum_{i=1}^n \mathbb{E}\Big[\rho\|x_{t+1} - x_t\|$$
$$\qquad\qquad + (1-\rho)(\|x_{t+1} - x_t\| + \|x_t - \alpha_t^i\|)\Big]$$
$$= \frac{1}{n} \sum_{i=1}^n \mathbb{E}\left[\|x_{t+1} - x_t\| + (1-\rho)\mathbb{E}\|x_t - \alpha_t^i\|\right] \qquad (14)$$

The first equality follows from the definition of $\alpha_{t+1}^i$ in Algorithm 4, while the inequality is just the triangle inequality. Using the above relationship and the bound in (12), we obtain

$$R_{t+1} \leq \mathbb{E}[F(x_t)] - \gamma \mathbb{E}[\mathcal{G}(x_t)] + \frac{LD^2\gamma^2}{2}$$
$$\qquad + \frac{LD\gamma\sqrt{n}}{\sqrt{b}} \frac{1}{n} \sum_{i=1}^n \mathbb{E}[\|x_t - \alpha_t^i\|] + c_{t+1}\mathbb{E}[\|x_{t+1} - x_t\|]$$
$$\qquad + c_{t+1}(1-\rho)\frac{1}{n} \sum_{i=1}^n \mathbb{E}[\|x_t - \alpha_t^i\|]$$
$$\leq R_t - \gamma \mathbb{E}[\mathcal{G}(x_t)] + \frac{LD^2\gamma^2}{2} + c_{t+1}D\gamma. \qquad (15)$$

The second inequality holds because: (a) $c_t = (1-\rho)c_{t+1} + (LD\gamma\sqrt{n})/\sqrt{b}$, and (b) $\|x_{t+1} - x_t\| = \gamma\|v_t - x_t\| \leq D\gamma$ (due to our bound on the diameter of the set $\Omega$). Telescoping over all the iterations, we see that

$$R_T \leq R_0 - \gamma \sum_{t=0}^{T-1} \mathbb{E}[\mathcal{G}(x_t)] + \frac{TLD^2\gamma^2}{2} + D\gamma \sum_{t=1}^{T} c_t$$

$$\leq R_0 - \gamma \sum_{t=0}^{T-1} \mathbb{E}[\mathcal{G}(x_t)] + \frac{TLD^2\gamma^2}{2} + \frac{LD^2\gamma^2\sqrt{n}}{\rho\sqrt{b}}$$

$$\leq R_0 - \gamma \sum_{t=0}^{T-1} \mathbb{E}[\mathcal{G}(x_t)] + \frac{TLD^2\gamma^2}{2} + \frac{2LD^2\gamma^2 n^{3/2}}{b^{3/2}}.$$

The second inequality follows form the fact that $\sum_{t=1}^{T} c_t \leq LD\gamma\sqrt{n}/(\rho\sqrt{b})$. This can, in turn, be obtained from the recursion $c_t = (1-\rho)c_{t+1} + (LD\gamma\sqrt{n})/\sqrt{b}$ and $c_T = 0$. The third inequality is due to the bound on $\rho$ in (13). Rearranging the above inequality and using the definition of $x_a$ from Algorithm 4, we finally obtain the bound

$$\mathbb{E}[\mathcal{G}(x_a)] \leq \frac{F(x_0) - \mathbb{E}[F(x_T)]}{T\gamma} + \frac{LD^2\gamma}{2} + \frac{2LD^2\gamma n^{3/2}}{Tb^{3/2}}$$

$$\leq \frac{F(x_0) - F(x^*)}{T\gamma} + LD^2\gamma\theta(b, n, T).$$

The first inequality uses the fact that $c_T = 0$ and $\alpha_0^i = x_0$ (in Algorithm 4). The second inequality uses the optimality of $x^*$ and the definition of $\theta(b, n, T)$. Using the setting of $\gamma$ in the theorem statement, we obtain the desired result. □

**Corollary 3.** *Assume $T \geq n$. Under the settings of Theorem 4 and with $b = \lceil n^{1/3} \rceil$, the IFO and LO complexity of Algorithm 2 are $O(n + n^{1/3}/\epsilon^2)$ and $O(1/\epsilon^2)$, respectively.*

*Proof.* First, observe that for $T \geq n$ and $b = \lceil n^{1/3} \rceil$, $\theta(b, n, T) \leq 5/2$ in Theorem 4. Thus, the IFO complexity is $O(n + n^{1/3}/\epsilon^2)$. Furthermore, since each iteration requires just $O(1)$ LO calls, the LO complexity is $O(1/\epsilon^2)$. □

Notably, the IFO complexity of SAGAFW is lower than that of SVFW. Moreover, if $T \geq n^{3/2}$ and $b = 1$, then we have $\theta(b, n, T) \leq 5/2$, in which case the IFO complexity is $O(n^{3/2} + 1/\epsilon^2)$.

## 4 Variance Reduction in Stochastic Setting

In this section, we improve the convergence rates in the stochastic setting using variance reduction techniques. The key idea is to first obtain samples $\{z_i\}$ are chosen independently according to the distribution $\mathcal{P}$ and then use SVFW or SAGAFW, described in this paper, on the finite-sum problem over these samples. The pseudocode for the SVFW and SAGAFW variants for stochastic setting (SVFW-S and SAGAFW-S respectively) are provided in Figure 2. The following is the key result regarding the convergence rates of SVFW-S and SAGAFW-S.

**Theorem 5.** *Consider the stochastic setting of (1) where $f$ is $G$-Lipschitz and $f(., z)$ is $L$-smooth for all $z \in \Xi$. Suppose $B = T$ and $\gamma = \sqrt{\frac{F(x_0) - F(x^*)}{TLD^2\beta}}$ (for SVFW-S and SAGAFW-S). Then the output of SVFW-S and SAGAFW-S satisfy the following:*

$$\mathbb{E}[\mathcal{G}(x_a)] \leq \frac{2D}{\sqrt{T\beta}}\sqrt{L(F(x_0) - F(x^*))(1+\beta)} + \frac{GD}{\sqrt{T}} \qquad (16)$$

| **SVFW-S:**$(x_0, T, B, \gamma)$ | **SAGAFW-S:**$(x_0, T, B, \gamma)$ |
|---|---|
| Randomly sample $z_1, \cdots, z_B \sim \mathcal{P}$ | Randomly sample $z_1, \cdots, z_B \sim \mathcal{P}$ |
| Let finite-sum $\hat{F}(x) = \frac{1}{B} \sum_{i=1}^{B} f(x, z_i)$ | Let finite-sum $\hat{F}(x) = \frac{1}{B} \sum_{i=1}^{B} f(x, z_i)$ |
| Output SVFW$(x_0, T, B^{1/3}, \gamma, \lceil B^{2/3} \rceil)$ applied on the function $\hat{F}$ | Output SAGAFW$(x_0, T, \gamma, \lceil 3B^{1/3} \rceil)$ applied on the function $\hat{F}$ |

Figure 2: SVFW-S and SAGAFW-S variants for the stochastic setting.

*Proof.* Consider the finite-sum $\hat{F}(x) = \frac{1}{B} \sum_{i=1}^{B} f(x, z_i)$ where $z_1, \cdots, z_B \sim \mathcal{P}$. We use the following notation:

$$\hat{\mathcal{G}}(x) = \max_{v \in \Omega} \langle v - x, -\nabla \hat{F}(x) \rangle.$$

Let $\bar{v}_t^{s+1} = \arg\max_{v \in \Omega} \langle v - x_t^{s+1}, -\nabla F(x_t^{s+1}) \rangle$ and $\hat{v}_t^{s+1} = \arg\max_{v \in \Omega} \langle v - x_t^{s+1}, -\nabla \hat{F}(x_t^{s+1}) \rangle$. We first observe the following key relationship for SVFW:

$$\mathbb{E}[\mathcal{G}(x_t^{s+1}) - \hat{\mathcal{G}}(x_t^{s+1})] = \mathbb{E}[\langle \bar{v}_t^{s+1} - x_t^{s+1}, -\nabla F(x_t^{s+1}) \rangle] - \mathbb{E}[\langle \hat{v}_t^{s+1} - x_t^{s+1}, -\nabla \hat{F}(x_t^{s+1}) \rangle]$$
$$\leq \mathbb{E}[\langle \bar{v}_t^{s+1} - x_t^{s+1}, -\nabla F(x_t^{s+1}) \rangle] - \mathbb{E}[\langle \bar{v}_t^{s+1} - x_t^{s+1}, -\nabla \hat{F}(x_t^{s+1}) \rangle]$$
$$\leq \mathbb{E}[\langle \bar{v}_t^{s+1} - x_t^{s+1}, \nabla \hat{F}(x_t^{s+1}) - \nabla F(x_t^{s+1}) \rangle]$$
$$\leq D\mathbb{E}[\|\nabla \hat{F}(x_t^{s+1}) - \nabla F(x_t^{s+1})\|] \leq \frac{GD}{\sqrt{T}}.$$

The first inequality is due to the optimality of $\hat{v}_t^{s+1}$. The third inequality follows from Cauchy-Schwarz inequality. The last inequality is due to Lemma 2. Adding the above inequality across all the iterations and epochs, we get:

$$\mathbb{E}[\mathcal{G}(x_a)] \leq \mathbb{E}[\hat{\mathcal{G}}(x_a)] + \frac{GD}{\sqrt{T}}.$$

Using the bound on $\mathbb{E}[\hat{\mathcal{G}}(x_a)]$ in Theorem 3 (here, recall we are using SVFW on $\hat{F}$) in the above inequality, we get the desired result. The proof for SAGAFW-S is similar. □

The following corollary on the complexity of SVFW-S and SAGAFW-S is immediate consequence of the above result.

**Corollary 4.** *Under the setting of Theorem 5, the SFO complexity of* SVFW-S *and* SAGAFW-S *(in Figure 2) are $O(1/\epsilon^{10/3})$ and $O(1/\epsilon^{8/3})$, respectively. The LO complexity of both* SVFW-S *and* SAGAFW-S *is $O(1/\epsilon^2)$.*

*Proof.* The proof follows from the fact that $B = T$, $b = \lceil B^{2/3} \rceil$ (in SVFW-S) and $b = \lceil 3B^{1/3} \rceil$ (in SAGAFW-S) and IFO complexities of SVFW and SAGAFW given in Corollary 2 and 3 respectively. □

By comparing Corollary 4 with Corollary 1, we see that SVFW-S and SAGAFW-S have better SFO complexity than SFW.

## 5 Discussion

It is important to remark on the complexity results derived in this paper. For the stochastic setting, we showed that the SFO and LO complexity of SFW are $O(1/\epsilon^4)$ and $O(1/\epsilon^2)$, respectively. At first glance, these rates might appear worse than those obtained for nonconvex SGD (see (Ghadimi and Lan, 2013)). However, it is important to note that the convergence criterion used in our paper is different from the one used in (Ghadimi and Lan, 2013). It is an important piece of future work to understand the precise relationship between these convergence criteria. Furthermore, the convergence rates in this paper are similar to those obtained for online Frank-Wolfe (Hazan and Kale, 2012) and only slightly worse than those obtained for stochastic Frank-Wolfe in the *convex setting* (Hazan and Luo, 2016). We, further, improved the convergence rate of SFW by using variance reduction ideas in the stochastic setting (SVFW-S and SAGAFW-S algorithms in Section 4). Understanding the tightness of these rates is an interesting open problem left as future work.

For the finite-sum setting, while the complexity results of SFW still hold, we obtained significantly faster convergence rates by using variance reduction techniques. The dependence of IFO and LO complexity of nonconvex SVFW and SAGAFW, on $\epsilon$ is $O(1/\epsilon^2)$, which matches the classical Frank-Wolfe algorithm (Lacoste-Julien, 2016). However, SVFW and SAGAFW exhibit a much weaker dependence on $n$ than FW; wherein, they are provably faster than the classical Frank-Wolfe by a factor of $n^{1/3}$ and $n^{2/3}$, respectively. Similar (but not same) benefits have also been reported for nonconvex SVRG and SAGA over gradient descent (Reddi et al., 2016a;b). Interestingly, there appears to be a gap between the convergence rates of SVFW and SAGAFW. Whether this gap is an artifact of our analysis or has deeper reasons remains open.

We conclude with a remark on a subtle point regarding the step size $\gamma$. The step size $\gamma$ in Theorems 2, 3, and 4 requires knowledge of parameters like $L$, $D$ and $F(x) - F(x^*)$. Typically, an estimate of the these values suffices in practice. In absence of such knowledge, one can completely eliminate this dependence of $\gamma$ on these parameters by simply choosing $\beta = \frac{2(F(x_0) - F(x^*))}{LD^2}$. Fortunately, this comes at the cost of only slightly worse *constants* in the convergence rate.

# Appendix

The following bound on the value of functions with Lipschitz continuous gradients is classical (see e.g., (Nesterov, 2003)).

**Lemma 1.** *If $f : \mathbb{R}^d \to \mathbb{R}$ is $L$-smooth, then*

$$f(x) \leq f(y) + \langle \nabla f(y), x - y \rangle + \frac{L}{2}\|x - y\|^2,$$

*for all $x, y \in \mathbb{R}^d$.*

The following lemma is useful for bounding the variance of the gradient estimate used in the stochastic setting.

**Lemma 2.** *Suppose the function $F(x) = \mathbb{E}_z[f(x,z)]$ where $z$ is a random variable with distribution $\mathcal{P}$ and support $\Xi$, and $\max_{z \in \Xi} \|\nabla f(x,z)\| \leq G$ for all $x \in \Omega$. Also, let $\bar{\nabla}_x = \frac{1}{b} \sum_{i \in I_t} \nabla f(x, z_i)$ where $\{z_i\}_{i=1}^b$ are i.i.d. samples from the distribution $\mathcal{P}$. Then, the following holds for any $x \in \Omega$:*

$$\mathbb{E}[\|\bar{\nabla}_x - \nabla F(x)\|] \leq \frac{G}{\sqrt{b}}.$$

*Proof.* The proof follows from a simple application of Lemma 5 and Jensen's inequality. □

The following result is useful for bounding the variance of the updates of SVFW and follows from a slight modification of a result in (Reddi et al., 2016a). We give the proof here for completeness.

**Lemma 3** (Reddi et al., 2016a). *Let $\tilde{\nabla}_t = \frac{1}{b_t}(\sum_{i \in I_t} \nabla f_i(x_t^{s+1}) - f_i(\tilde{x}^s) + \tilde{g}^s)$ in Algorithm 3. For the iterates $x_t^{s+1}$ and $\tilde{x}^s$ where $t \in \{0, \ldots, m-1\}$ and $s \in \{0, \ldots, S-1\}$ in Algorithm 3, the following inequality holds:*

$$\mathbb{E}_{I_t}[\|\nabla F(x_t^{s+1}) - \tilde{\nabla}_t\|] \leq \frac{L}{\sqrt{b_t}} \|x_t^{s+1} - \tilde{x}^s\|.$$

*Proof.* For the ease of exposition, we first define

$$\zeta_t^{s+1} = \frac{1}{|I_t|} \sum_{i \in I_t} \left(\nabla f_i(x_t^{s+1}) - \nabla f_i(\tilde{x}^s)\right).$$

Using this notation, we then obtain the following:

$$\mathbb{E}_{I_t}[\|\nabla F(x_t^{s+1}) - \tilde{\nabla}_t\|^2]$$
$$= \mathbb{E}_{I_t}[\|\zeta_t^{s+1} + \nabla F(\tilde{x}^s) - \nabla F(x_t^{s+1})\|^2]$$
$$= \mathbb{E}_{I_t}[\|\zeta_t^{s+1} - \mathbb{E}_{I_t}[\zeta_t^{s+1}]\|^2]$$
$$= \frac{1}{b_t^2} \mathbb{E}_{I_t}\left[\left\|\sum_{i \in I_t} \left(\nabla f_i(x_t^{s+1}) - \nabla f_i(\tilde{x}^s) - \mathbb{E}_{I_t}[\zeta_t^{s+1}]\right)\right\|^2\right].$$

The second equality is due to the fact that $\mathbb{E}_{I_t}[\zeta_t^{s+1}] = \nabla F(x_t^{s+1}) - \nabla F(\tilde{x}^s)$. From the above relationship, we get

$$\mathbb{E}_{I_t}[\|\nabla F(x_t^{s+1}) - \tilde{\nabla}_t\|^2]$$
$$\leq \frac{1}{b_t} \mathbb{E}_{I_t}\left[\sum_{i \in I_t} \|\nabla f_i(x_t^{s+1}) - \nabla f_i(\tilde{x}^s) - \mathbb{E}_{I_t}[\zeta_t^{s+1}]\|^2\right]$$
$$\leq \frac{1}{b_t} \mathbb{E}_{I_t}\left[\sum_{i \in I_t} \|\nabla f_i(x_t^{s+1}) - \nabla f_i(\tilde{x}^s)\|^2\right]$$
$$\leq \frac{L^2}{b_t} \|x_t^{s+1} - \tilde{x}^s\|^2.$$

The first inequality follows from Lemma 5. The second inequality is due to the fact that for a random variable $\zeta$, $\mathbb{E}[\|\zeta - \mathbb{E}[\zeta]\|^2] \leq \mathbb{E}[\|\zeta\|^2]$. The last inequality follows from $L$-smoothness of $f_i$. The result follows from a simple application of Jensen's inequality to the inequality above. □

The following result is important for bounding the variance in SAGAFW. The key difference from Lemma 3 is that the variance term in SAGAFW involves $\alpha_t^i$. Again, we provide the proof for completeness.

**Lemma 4.** *Let $\check{\nabla}_t = \frac{1}{b_t}(\sum_{i \in I_t} \nabla f_i(x_t) - f_i(\alpha_t^i) + g_t)$ in Algorithm 4. For the iterates $x_t, v_t$ and $\{\alpha_t^i\}_{i=1}^n$ where $t \in \{0, \ldots, T-1\}$ in Algorithm 4, we have the inequality*

$$\mathbb{E}_{I_t}[\|\nabla F(x_t) - \check{\nabla}_t\|] \leq \frac{L}{\sqrt{b_t}} \sum_{i=1}^n \frac{1}{\sqrt{n}} \|x_t - \alpha_t^i\|.$$

*Proof.* As before we first define the quantity

$$\zeta_t = \frac{1}{|I_t|} \sum_{i \in I_t} \left( \nabla f_i(x_t) - \nabla f_i(\alpha_t^i) \right).$$

With this notation, we then obtain the equality

$$\mathbb{E}_{I_t}[\|\nabla F(x_t) - \check{\nabla}_t\|^2]$$
$$= \mathbb{E}_{I_t}\left[\left\|\zeta_t + \frac{1}{n}\sum_{i=1}^n \nabla f_i(\alpha_t^i) - \nabla F(x_t)\right\|^2\right] = \mathbb{E}_{I_t}[\|\zeta_t - \mathbb{E}_{I_t}[\zeta_t]\|^2]$$
$$= \frac{1}{b^2}\mathbb{E}_{I_t}\left[\left\|\sum_{i \in I_t}\left(\nabla f_i(x_t) - \nabla f_i(\alpha_t^i) - \mathbb{E}_{I_t}[\zeta_t]\right)\right\|^2\right].$$

The second equality follows from the fact that $\mathbb{E}_{I_t}[\zeta_t] = \nabla F(x_t) - \frac{1}{n}\sum_{i=1}^n \nabla f_i(\alpha_t^i)$. From the above inequality, we get the following bound:

$$\mathbb{E}_{I_t}[\|\nabla F(x_t) - \check{\nabla}_t\|^2]$$
$$\leq \frac{1}{b_t}\mathbb{E}_{I_t}\left[\sum_{i \in I_t} \|\nabla f_i(x_t) - \nabla f_i(\alpha_t^i) - \mathbb{E}_{I_t}[\zeta_t]\|^2\right]$$
$$\leq \frac{1}{b_t}\mathbb{E}_{I_t}\left[\sum_{i \in I_t} \|\nabla f_i(x_t) - \nabla f_i(\alpha_t^i)\|^2\right] \leq \frac{L^2}{nb_t}\sum_{i=1}^n \|x_t - \alpha_t^i\|^2.$$

The first inequality is due to Lemma 5, while the second inequality holds because for a random variable $\zeta$, $\mathbb{E}[\|\zeta - \mathbb{E}[\zeta]\|^2] \leq \mathbb{E}[\|\zeta\|^2]$. The last inequality is from $L$-smoothness of $f_i$ ($i \in [n]$) and uniform randomness of the set $I_t$. By applying Jensen's inequality, we get the desired result. □

**Lemma 5.** *For random variables $z_1, \ldots, z_r$ that are independent and have mean 0, we have*

$$\mathbb{E}\left[\|z_1 + \ldots + z_r\|^2\right] = \mathbb{E}\left[\|z_1\|^2 + \ldots + \|z_r\|^2\right].$$

*Proof.* Expanding the left hand side we have

$$\mathbb{E}\left[\|z_1 + \ldots + z_r\|^2\right] = \sum_{i,j=1}^r \mathbb{E}[z_i z_j] = \mathbb{E}\left[\sum_{i=1}^r \|z_i\|^2\right];$$

the second equality here follows from the our hypothesis. □